\documentclass[a4paper,12pt]{amsart}
\usepackage{amssymb,amstext,amsmath,amscd,amsthm,
amsfonts,enumerate,graphicx,latexsym}
\usepackage{url}
\usepackage[usenames]{color}
\definecolor{gr}{rgb}{0.1, .5 , .10}
\usepackage[all]{xy}

\newtheorem{theorem*}{Theorem}
\newtheorem*{theorem**}{Theorem}

\newtheorem{corollary*}[theorem*]{Corollary}
\newtheorem*{corollary**}{Corollary}

\theoremstyle{definition}

\newtheorem*{question*}{Question}

\newtheorem*{conjecture*}{Conjecture}

\newtheorem*{example**}{Example}

\newtheorem*{notation*}{Notation}

\newtheorem*{claim*}{Claim}

\numberwithin{equation}{theorem}

\def\Hom{\operatorname{Hom}}

\def\End{\operatorname{End}}

\def\thick{\operatorname{\mathsf{thick}}}

\def\D{\operatorname{\mathsf{D}}}

\def\proj{\operatorname{\mathsf{proj}}}

\def\Kb{\mathsf{K^b}}

\def\per{\operatorname{\mathsf{per}}}
\def\Tria{\operatorname{\mathsf{Tria}}}

\def\A{\mathcal{A}}
\def\B{\mathcal{B}}
\def\C{\mathcal{C}}

\def\DD{\mathcal{D}}
\def\T{\mathcal{T}}

\def\X{\mathcal{X}}
\def\Y{\mathcal{Y}}

\def\E{\mathcal{E}}

\def\KK{\mathcal{K}}

\begin{document}
\setlength{\baselineskip}{15pt}
\title[A note on a homological epimorphism]{A note on a homological epimorphism whose source is a piecewise hereditary algebra}
\author{Takuma Aihara}
\address{Department of Mathematics, Tokyo Gakugei University, 4-1-1 Nukuikita-machi, Koganei, Tokyo 184-8501, Japan}
\email{aihara@u-gakugei.ac.jp}

\thanks{TA was partly supported by JSPS KAKENHI Grant Number JP25K06950.}
\begin{abstract}
We study a homological epimorphism of proper dg algebras and show that if the source is a finite dimensional piecewise hereditary algebra, then the target is derived equivalent to some piecewise hereditary algebra.
\end{abstract}
\maketitle
\section{Introduction and results}

Let $\Lambda$ be a finite dimensional algebra over an algebraically closed field $K$ and $S$ a (basic) presilting object in the perfect derived category $\Kb(\proj\Lambda)$ of $\Lambda$.
Then there exists a recollement $(\D(\A), \D(\Lambda), \Tria S)$, where $\D(-)$ stands for the (unbounded) derived category and $\A$ is the dg endomorphism algebra of $\Lambda$ in a dg enhancement of $\Kb(\proj\Lambda)/\thick S$ \cite{AH, KY}.
A problem is what the dg algebra $\A$ is.
We know that for each $i$, the cohomology $H^i(\A)$ is finite dimensional over $K$.
For example, if $S=e\Lambda$ for an idempotent $e$ of $\Lambda$, then $\A$ is quasi-isomorphic to $\widetilde{\Lambda}/\widetilde{\Lambda}\widetilde{e}\widetilde{\Lambda}$ \cite{KY}.
Here, $\widetilde{\Lambda}$ is a dg path algebra quasi-isomorphic to $\Lambda$ \cite{O}.
However, it is quite difficult to answer the question in general.

Assuming $\Lambda$ to be piecewise hereditary, we obtain the main theorem of this paper.

\begin{theorem**}
Keeping above, $\A$ is derived equivalent to some piecewise hereditary algebra.
Moreover, if $\Lambda$ is of Dynkin type, then so is the algebra.
\end{theorem**}

\begin{example**}
Let $\varphi: \Lambda\to \B$ be a homological epimorphism of dg algebras with $\B$ proper.
Then, we have a recollement $(\D(\B), \D(\Lambda), \X)$ \cite{NS}.
Now, let us suppose that $\Lambda$ is piecewise hereditary.
Since the dg algebra $\B$ is compact in $\D(\Lambda)$, it is seen that $\X$ is compactly generated.
By \cite{AHMV}, we obtain a presilting object $S$ of $\Kb(\proj\Lambda)$ with $\X=\Tria S$.
Hence, the theorem leads to the fact that $\B$ is derived equivalent to a piecewise hereditary algebra $\Gamma$; if in addition, $\Lambda$ is of Dynkin type, then so is $\Gamma$.
\end{example**}

By \cite[Corollary 2.17]{AH}, we recover a result of Xu--Yang on Bongartz-type lemma \cite{XY}.

\begin{corollary**}
Any presilting of a piecewise hereditary algebra can be completed to silting.
\end{corollary**}

\section{Proof of the theorem}

We may assume that $\Lambda$ is a hereditary algebra; when it is a canonical algebra, we regard it as a hereditary abelian category with tilting object.
Put $\T:=\Kb(\proj\Lambda)$.

Let $S:=S_1\oplus \cdots\oplus S_r$ be the indecomposable decomposition of the presilting object.
By \cite[Proposition 3.11]{AI}, under a suitable relabelling we obtain an exceptional sequence $(S_1, \cdots, S_r)$; that is, $S_i$ is pretilting-brick and $\Hom_\T(S_i, S_j[t])=0$ for $1\leq j<i\leq r$ and all integers $t$.
Then, there exists a semiorthogonal decomposition $\T=\DD *\C$ with $\C:=\thick S$ and $\DD:={}^\perp\C$.
Since we have an exact sequence $0\to \C\to\T\to \per(\A)\to0$ of triangulated categories (see \cite{AH} for example), one gets a triangle equivalence $\per(\A)\simeq \DD$.

By the induction on $r$, we show that $\DD$ admits a tilting object $M$ with $\End_\T(M)$ piecewise hereditary.
If $r=1$, then there is a (classical) tilting object $S\oplus M$ with $M\in\DD$  because of the classical Bongartz's lemma;
see also \cite[Proposition 2.16]{A}.
As $M$ is partial tilting, we see that $\End_\T(M)$ is piecewise hereditary \cite[Corollary 5.10]{AHKL}.

Let $r>1$ and $S':=S_1\oplus\cdots\oplus S_{r-1}$.
As $(S_1,\cdots,S_{r-1})$ is an exceptional sequence, we have a semiorthogonal decomposition $\T=\DD'*\C'$ with $\C':=\thick S'$ and $\DD':={}^\perp\C'$.
By the induction hypothesis, we observe that $\DD'$ possesses a tilting object $M'$ with $\End_\T(M')$ piecewise hereditary.
Since $S_r$ is a member of $\DD'$ and exceptional in $\DD'$, there exists a semiorthogonal decomposition $\DD'={}^{\perp_{\DD'}}\E*\E$ with $\E=\thick_{\DD'}S_r$ and ${}^{\perp_{\DD'}}\E$ admits a tilting object $M$ with $\End_\T(M)$ piecewise hereditary.
It is not difficult to check that ${}^{\perp_{\DD'}}\E=\DD\ (:={}^\perp(\thick S))$, whence we are done.

Finally, assume that $\Lambda$ is of Dynkin type.
Since $\Kb(\proj\Lambda)$ is silting-discrete, we derive from \cite[Theorem 2.15]{AH} that $\per(\A)\simeq \Kb(\proj\End_\T(M))$ is also silting-discrete.
This implies that $\End_\T(M)$ is of Dynkin type.

\section{Remarks}

\begin{enumerate}
\item The theorem is a derived category version of Crawly-Boevey's theorem \cite{CB}, and is a dg/silting version of Angeleri-Hugel--Koenig--Liu's theorem \cite{AHKL}.



\item If $\Lambda$ is a piecewise hereditary algebra, then the thick closure of $S$ is also triangle equivalent to the perfect derived category of some piecewise hereditary algebra.
In particular, if a recollement $(\Y, \D(\Lambda), \X)$ with $\X$ compactly generated exists, then it is that of piecewise hereditary algebras; i.e., $(\Y, \D(\Lambda), \X)\simeq (\D(\Gamma), \D(\Lambda), \D(\Sigma))$.

\item We can never drop the assumption of $\Lambda$ being a piecewise hereditary algebra.
In fact, it is known that the dg algebra $\A$ does not necessarily admit a silting object \cite{CJS, Kr, Ka}; all the algebras $\Lambda$ given in the papers have finite global dimension.
Thus, even if the global dimension of $\Lambda$ is finite, $\A$ is not always derived equivalent to a finite dimensional/nonpositive (connective) dg algebra.

\item The center in a recollement does not necessarily inherit the heredity property from the both sides.
We close this paper by giving an example below.
\end{enumerate}

\begin{example**}
Let $\KK$ be the dg algebra presented by the quiver $\xymatrix{1 \ar@<2pt>[r]^x\ar@{-->}@<-2pt>[r]_y & 2}$ with $\deg(x)=0$ and $\deg(y)=\ell$ (trivial differential). 
Taking the primitive idempotent $e$ of $\KK$ corresponding to the vertex 2, we have a recollement $(\D(\KK/\KK e\KK), \D(\KK), \Tria(e\KK))$.
One observes that $\D(\KK/\KK e\KK)\simeq \D(K)\simeq \Tria(e\KK)$; so, the both sides are hereditary.
However, the center is not derived equivalent to an \emph{ordinary} algebra unless $|\ell|\leq 1$ \cite[Corollary 3.2]{LY}.

If $|\ell|=1$, then by \cite[Corollary 3.10]{LY} we see that $\KK$ is derived equivalent to the algebra given by the quiver $\xymatrix{1 \ar@<2pt>[r]^y & 2 \ar@<2pt>[l]^z}$ with relation $zy=0$; this is not piecewise hereditary, either.
%

Applying the Koszul dual, let us suppose $\ell\leq0$.
Then, we observe that $\KK$ has global dimension $1-\ell$ in the sense of Minamoto \cite{M}.
\end{example**}


\end{document}